\documentclass[12pt]{article}
\usepackage{oldlfont}
\usepackage{enumerate}
\usepackage{hyperref}
\usepackage{amssymb}
\usepackage[left=25mm,right=25mm]{geometry}

\begin{document}
%
%

\pagestyle{headings}                                                    
\flushbottom                                                            

\makeatletter                                                           

\def\section{\@startsection {section}{1}{\z@}{-1.5ex plus -.5ex         
minus -.2ex}{1ex plus .2ex}{\large\bf}}                                 

\renewcommand{\thesection}{\arabic{section}.}                           
\renewcommand{\thesubsection}{\thesection\arabic{subsection}.}          
\renewcommand{\thesubsubsection}{\thesubsection\arabic{subsubsection}.} 
\renewcommand{\theparagraph}{\thesubsubsection\arabic{paragraph}.}      
\renewcommand{\thesubparagraph}{\theparagraph\arabic{subparagraph}.}    

\def\@thmcountersep{}                                                   

\long\def\@makecaption#1#2{\vskip 10pt \setbox\@tempboxa\hbox{#1. #2}   
   \ifdim \wd\@tempboxa >\hsize   
       #1. #2                  
     \else                        
       \hbox to\hsize{\hfil\box\@tempb oxa\hfil}                         
   \fi}                                                                 

\def\ps@headings{                                                       
 \def\@oddhead{\footnotesize\rm\hfill\runninghead\hfill}               
 \def\@evenhead{\@oddhead}                                              
 \def\@oddfoot{\rm\hfill\thepage\hfill}\def\@evenfoot{\@oddfoot} }      

\newcommand{\FF}{{\rm I\kern-1.6pt {\rm F}}}
\newtheorem{thm}{Theorem}[section]
\newtheorem{lemma}[thm]{Lemma}
\newtheorem{cor}[thm]{Corollary}
\newtheorem{prop}[thm]{Proposition}
\newtheorem{num}[thm]{}
\newcommand{\pf}{\noindent {\em Proof.} \ }
\newcommand{\eop}{$_{\Box}$ \vspace{5mm} \relax}
\newcommand{\MSC}{{\rm Mathematics Subject Classification: }}
\newcommand{\Aut}{{\rm Aut}}
\newcommand{\PSL}{{\rm PSL}}
\newcommand{\PGL}{{\rm PGL}}
\newcommand{\SL}{{\rm SL}}
\newcommand{\GL}{{\rm GL}}
\newcommand{\AGL}{{\rm AGL}}
\newcommand{\ASL}{{\rm ASL}}
\newcommand{\SU}{{\rm SU}}
\newcommand{\PSU}{{\rm PSU}}
\newcommand{\Sym}{{\rm Sym}}
\newcommand{\nl}{\trianglelefteq}
\newcommand{\sdp}{\rtimes}
\renewcommand{\o}{\overline}

\providecommand{\keywords}[1]{\textbf{Keywords:} #1}
\newcommand\blfootnote[1]{%
  \begingroup
  \renewcommand\thefootnote{}\footnotetext{#1}%
  \addtocounter{footnote}{-1}%
  \endgroup}
\newcommand{\Addresses}{{
  \bigskip
  \footnotesize
\noindent
  \textsc{Dipartimento di Matematica e Informatica, Universit\`a della Calabria,
    87036 Arcavacata di Rende, Italy}\par\nopagebreak\noindent
  \textit{E-mail address}: \texttt{jozef.vanbon@unical.it}
}}

\title{\bf A locally $5$-arc transitive graph related to $\PSU_3(8)$}

\author{J. van Bon}

\date{\ }

\maketitle

\Addresses

\def\runninghead{{\sc Locally $5$-arc transitive graphs}}
\pagestyle{headings}

\begin{abstract}
We construct a locally 5-arc transitive $G$-graph $\Delta$, where
$G \in \{\PSU_3(8) \sdp C_2, \PSU_3(8) \sdp C_6\}$.
The corresponding vertex stabilizer amalgams are of local
characteristic $3$ but are not weak $(B,N)$-pairs. They are the 
first examples of this kind where there exists a vertex $z$ such that
the extension of $G_z/O_p(G_z^{[1]})$ over $O_p(G_z^{[1]})$ is non-split.
\end{abstract}


\par\smallskip

\section{Introduction.}

A $G$-graph is a simple graph $\Delta$ together with a subgroup $G \leq \Aut(\Delta)$.
An $s$-arc in $\Delta$ is a path $(x_0, x_1, \dots ,x_s)$ with
$x_{i-1}\neq x_{i+1}$ for $1\leq i \leq s-1$. 
A $G$-graph $\Delta$ is called {\em locally $s$-arc transitive} if for any
vertex $z$ the stabilizer $G_z$ is transitive on the set of $s$-arcs
originating at $z$. The study of locally $s$-arc transitive graphs was initiated by
Tutte \cite{T1, T2} who showed that if $\Delta$ is a 
finite trivalent $G$-graph for which $G$ is
transitive on the vertex set of $\Delta$, then $s \leq 5$.
Assume from now on that $\Delta$ is a locally $s$-arc transitive $G$-graph with 
$s\geq 1$, and that the valency of each vertex is at least 3.
Since $s \geq 1$ the group 
$G$ is transitive on the set of edges $E \Delta$ and thus
has at most two orbits on the vertex set $V \Delta$. 
Let $\{x_1,x_2\}$ be an edge. The triple $(G_{x_1},G_{x_2};G_{x_1,x_2})$ is
called the {\em vertex stabilizer amalgam} with respect to the 1-arc $(x_1,x_2)$. The $G$-graph will be called 
{\em locally finite} if the groups $G_{x_1}$ and $G_{x_2}$ are both finite groups.
Tutte's result was generalized by Goldschmidt \cite{Gold} who showed that if
$\Delta$ is a trivalent locally finite $G$-graph,
then $s\leq 7$ and the vertex stabilizer amalgam belongs to one of 15 isomorphism types.
There is no hope to classify the corresponding graphs though. Likewise for an arbitrary 
locally finite and locally $s$-arc transitive $G$-graph
the best one can hope for
is to bound $s$ and determine the possible isomorphism type of the vertex stabilizer amalgam for certain values of $s$.
The first problem was solved in \cite{bonstell} where the bound $s \leq 9$ was obtained.
Vertex stabilizer amalgams of locally finite and locally $s$-arc transitive $G$-graphs can have many different isomorphism types
and a classification of these amalgams for all possible values of $s$ is impossible.
However, for $s \geq 4$ the structure of the vertex stabilizer amalgam seems to be highly restricted. 

Before continuing our discussion
we need to introduce some more notation and definitions.
Let $z \in V\Delta$.
The kernel of action of $G_z$ on the set of neighbors $\Delta (z)$ 
will be denoted by $G_z^{[1]}$.
A locally finite $G$-graph $\Delta$ is called of
\begin{enumerate}[$\cdot$]
\item {\em local characteristic $p$} with respect to $G$,
if there exists a prime $p$ such that
$$C_{G_{x_i}}(O_p(G_{z}^{[1]})) \leq O_p(G_{z}^{[1]}), \hbox{ for all } z\in V\Delta;$$
\item {\em pushing up type} with respect to the 1-arc $(x_1,x_2)$ and the prime $p$, if $\Delta$ is of local characteristic $p$ 
respect to $G$ and $O_p(G_{x_1}^{[1]}) \leq O_p(G_{x_2}^{[1]})$.
\end{enumerate}

The most general result about the structure of vertex stabilizer of
locally finite $G$-graphs
is due to Thompson \cite{Thom} and Wielandt \cite{Wiel}. 
This result, known as the Thompson-Wielandt Theorem, 
states that -- under some mild conditions which are satisfied if $s \geq 2$ -- here exists a prime $p$ such that 
$G_{x_1,x_2}^{[1]}$ is a $p$-group or $G_{x_1,x_2}^{[1]}=G_{x_1}^{[2]}$ and  $G_{x_2}^{[2]}=G_{x_2}^{[3]}$ is a $p$-group
(see \cite{TW}). It
 was shown in \cite{kern} that $G_{x_1,x_2}^{[1]}$ is non-trivial when $s \geq 4$. 

The structure of the vertex stabilizer amalgam of $G$-graphs $\Delta$ with $s \geq 4$ appearing in the literature can be divided into three classes:

\begin{enumerate}[$\cdot$]
\item $\Delta$ is not of local characteristic $p$ and $s=5$;
\item $\Delta$ is of pushing up type and  $s=5$;
\item the vertex stabilizer amalgam of $\Delta$ is a weak $(B,N)$-pair.
\end{enumerate}
(see \cite{DS} for the definition of a weak $(B,N)$-pair and their classification.)
For locally finite $G$-graphs 
with $s \geq 4$ that are not of local characteristic $p$, the possible vertex stabilizer amalgams were determined in \cite{NonLocalChar}.
It is conjectured that the vertex stabilizer amalgam of a $G$-graph with $s \geq 4$ which
is of local characteristic $p$ but not of pushing up type is a weak $(B,N)$-pair.
This was shown to be true under the stronger assumption of $s \geq 6$ in \cite{bonstell}.

This paper forms part of the project \cite{PvB,PushingUpType1,PushingUpType} of determining the isomorphism types of vertex stabilizer amalgams of locally $s$-arc transitive graphs of pushing up type with $s \geq 4$. 
Here we construct a graph $\Delta$ and show that there are
subgroups $H \leq K \leq \Aut(\Delta)$ such that, for $G$ equal to $H$ and $K$,
$\Delta$ is a locally finite $G$-graph of pushing up type with respect to
the 1-arc $\{x_1,x_2\}$ and the prime 3.
The corresponding vertex stabilizer amalgams are of shape ${\cal D}_2$ and ${\cal E}_2$, respectively (see section 2 for their definitions).
In particular this shows that these two amalgams, which appeared in the conclusion of \cite[1.1]{PushingUpType1}, do occur as
vertex stabilizer amalgams of a $G$-graph. The graph $\Delta$ is the first example of a locally 5-arc transitive $G$-graph with a vertex stabilizer amalgam of this kind.

\begin{thm}\label{main1}
Let $\Delta$, $H$ and $K$ be the graph and the two groups constructed in Section 3,
respectively.
Then $H \leq K = \Aut(\Delta)$. Moreover, for $G\in\{H,K\}$,  $\Delta$ is a locally 5-arc transitive $G$-graph  
of local characteristic $3$ such that for an edge $\{x_1,x_2\}$ the following hold:
\begin{enumerate}[(i)]
\item $|\Delta(x_1)|=4$ and $|\Delta(x_2)|=3$;
\item $\Delta$ is of pushing up type with respect to the 1-arc $(x_1,x_2)$ and prime 3;
\item $G_{x_2}$ is transitive on 6-arcs originating at $x_2$;
\item The vertex stabilizer amalgam ${\mathfrak A}=(H_{x_1},H_{x_2};H_{x_1,x_2})$ is of shape ${\cal D}_2$;
\item The vertex stabilizer amalgam ${\mathfrak B}=(K_{x_1},K_{x_2};K_{x_1,x_2})$ is of shape ${\cal E}_2$.
\end{enumerate}
Moreover, $H \cong \PSU_3(8) \sdp \langle \sigma^3 \rangle$ and $K \cong \PSU_3(8) \sdp \langle \sigma \rangle$, where $\sigma$ is a field automorphism
of order 6 of $\PSU_3(8)$.
\end{thm}

\par\medskip
For any positive integer $i$ we define
$$G_z^{[i]}=\{g \in G_z \mid x^g=x \hbox{ for all } x\in V\Delta \hbox{ with } d(z,x) \leq i\},$$
that is,  $G_z^{[i]}$ is the kernel of action $G_z$ on the set of vertices at distance 
at most $i$ from $z$.  The group induced by $G_z$ on $\Delta (z)$ will be denoted
by $G_z^{\Delta (z)}$, thus $G_z^{\Delta (z)} \cong G_z /G_z^{[1]}$.
The following theorem gives the structure of the vertex stabilizers in
terms of the various kernels of action.

\begin{thm}\label{main2}
Let $H$, $K$ and $\Delta$ be as in Theorem 1.1.
Let $W_i=O_p(H_{x_i}^{[1]})$ and $\hat W_i = O_p(K_{x_i}^{[1]})$, $i \in \{1,2\}$.
Then $W_1$ is a elementary abelian group of order 9, 
$W_2$ is a special group of order 27 and exponent 3,
$\hat W_1 \cong W_1 \times C_3$ and $\hat W_2 \cong W_2 \times C_3$.
Moreover the following hold:
\begin{enumerate}[(i)]
\item $H_{x_1}^{\Delta (x_1)} \cong \Sym(4)$, $H_{x_1}^{[1]} \cong W_1 \sdp C_2$, $H_{x_1}^{[2]} \cong W_1$, $H_{x_1}^{[3]}=1$;
\item $H_{x_2}^{\Delta (x_2)} \cong \Sym(3)$, $H_{x_2}^{[1]} \cong W_2 \sdp C_2$, $H_{x_2}^{[2]}=H_{x_2}^{[3]} =Z(W_2)\cong C_3$, $H_{x_2}^{[4]}=1$;
\item $K_{x_1}^{\Delta (x_1)} \cong \Sym(4)$, $K_{x_1}^{[1]} \cong \hat W_1 \sdp C_2$, $K_{x_1}^{[2]} \cong \hat W_1$, $K_{x_1}^{[3]} =K_{x_1}^{[4]}=Z(K_{x_1}) \cong C_3$, $K_{x_1}^{[5]}=1$;
\item $K_{x_2}^{\Delta (x_2)} \cong \Sym(3)$, $K_{x_2}^{[1]} \cong \hat W_2 \sdp C_2$, $K_{x_2}^{[2]}=K_{x_2}^{[3]} =Z(\hat W_2) \cong C_3 \times C_3$, $K_{x_2}^{[4]}=1$.
\end{enumerate}
\end{thm}

\par\smallskip\noindent
{\bf Remark}. $K$-graph $\Delta$ is the first example of a graph of pushing up type with 
$K_{x_1}^{[3]}\neq 1$.

\par\smallskip
The organization of the paper is as follows. In section 2 we define the two amalgams
occurring in the conclusion of \ref{main1}. Section 3 is devoted to the proof of the two main theorems.

\par\smallskip
\section{Notation and definitions}

Our group theoretic notation is as follows.
Let $\ASL_2(3)$ and $\AGL_2(3)$, denote the semi-direct product of $\SL_2(3)$ and $\GL_2(3)$, respectively,
with its natural 2-dimensional module.
Let $\AGL_1(3)$ denote the semi-direct product of $\GL_1(3)$ with its natural 1-dimensional module.

With $SP_2$ we denote the extraspecial group of order 27 and of exponent 9.
Note that $\PGL_2(3)\cong \Sym(4)$ and $\AGL_1(3)\cong \Sym(3)$.

Let $S \in Syl_3(\AGL_2(3))$, $V=O_3(\AGL_2(3))$ and $V_0=C_V(S)=Z(S)$.
Let $\AGL_2(3,S) = N_{\AGL_2(3)}(S)\cong E_9 \sdp (C_2 \times \AGL_1(3))$.
With $\AGL_2(3,S)^{\#}$ we denote the unique subgroup $X$ of index 2 in $\AGL_2(3,S)$ with
$$X=C_X(V_0)C_X(V/V_0), \, C_X(V/V_0)/S = 1 \hbox{ and } C_X(V_0)/S\cong \GL_1(3).$$
With $\AGL_2(3,S)^*$ we denote the unique subgroup $Y$ of index 2 in $\AGL_2(3,S)$ with
$$Y=C_Y(V_0)C_Y(V/V_0), \, C_Y(V/V_0)/S \cong \GL_1(3) \hbox{ and } C_Y(V_0)/S=1.$$

For a brief outline of the theory of amalgams see \cite{NonLocalChar} whose notation and set up we follow.

\par\noindent
An amalgam ${\mathfrak A}: H_1 \leftarrow H_{1,2} \rightarrow H_2$ is called
{\em finite} if $H_1$ and $H_2$ are finite groups

\par\smallskip
Let ${\mathfrak A}: H_1 \leftarrow H_{1,2} \rightarrow H_2$ be a finite amalgam.
Following \cite{NonLocalChar}
we define $T_i$ to be the largest subgroup of $H_{1,2}$ such that $T_i \nl H_i$.
Let $X$ be the unique inclusion-minimal subgroup of $H_{1,2}$ satisfying
$$T_1T_2 \leq X \hbox{ and } X=H_{1,2}\cap \langle X^{H_i} \rangle, \hbox{ for } i=1,2.$$

\par\smallskip\noindent
We are now ready to define the two amalgams appearing in the conclusion of \ref{main1}.
These amalgams were first introduced in \cite{PushingUpType1}.

\par\smallskip\noindent
{\bf Definition} ({\em Amalgams of shape $\ASL_2(3,S)$}.)
Let ${\mathfrak A}: H_1 \leftarrow H_{1,2} \rightarrow H_2$ be a finite amalgam.
Then ${\mathfrak A}$ has shape $\AGL_2(3,S)$ if the following hold:
\begin{enumerate}[(i)]
\item $H_{1,2}=X$;
\item $T_2=C_X (Z(O_3(X)))$ and $H_2/Z(O_3(X)) \cong \AGL_2(3,S)$.
\end{enumerate}

\par\smallskip\noindent
{\bf Definition}
({\em Amalgams of shape ${\cal D}_2$ and ${\cal E}_2$}.)
Let ${\mathfrak A}$ be an amalgam of shape $\AGL_2(3,S)$.

\par\smallskip\noindent
We say that ${\mathfrak A}$ is of shape   ${\cal D}_2$
if the following hold:
\begin{enumerate}
\item $H_1 \cong \AGL_2(3)$ and $H_{1,2}  \cong \AGL_2(3,S)$;
\item $H_2 \cong O_3(X) \cdot (C_2 \times \AGL_1(3))$ and $T_2 \cong \AGL_2(3,S)^{\#}$;
\item $C_{O_3(H_2)}(T) \cong  C_9$, where $T \in Syl_2(T_2)$.
\end{enumerate}

\par\smallskip\noindent
We say that ${\mathfrak A}$ is of shape ${\cal E}_2$
if the following hold:
\begin{enumerate}
\item $H_1 \cong C_3 \times \AGL_2(3)$, $H_{1,2} \cong C_3 \times \AGL_2(3,S)$;
\item $H_2 \cong O_3(X) \cdot (C_2 \times \AGL_1(3))$ and $T_2 \cong C_3 \times \AGL_2(3,S)^{\#}$;
\item $Z(O_3(X)) \sdp H_2/C_{H_2}(Z(O_3(X))) \cong \AGL_2(3,S)^*$;
\item $C_{O_3(H_2)}(T) \cong  SP_2$, where $T \in Syl_2(T_2)$.
\end{enumerate}

\par\smallskip
\section{The construction of $\Delta$ and the proof of the main results}

In this section we first construct the graph $\Delta$
using some subgroups $K_1$, $K_2$ of $\Aut (\PSU_3(8))$ and then proceed proving the 
two main results of this paper.
The subgroups $K_1$ and $K_2$ can already be found in \cite{atlas}.
However the description there does not exhibit clearly
the structure of $O_3(K_1)$ and $O_3(K_2)$ which is needed show that the 
$G$-graph  $\Delta$ has $s=5$.
We therefor have included the details of the proofs.

Let $\zeta$ be a generator of $\FF_{64}^*$ and 
$\rho : \FF_{64} \rightarrow \FF_{64}$ be the standard Frobenius map 
$x \mapsto x^2$. Let $\tau=\rho^3$, so $\tau$ has order $2$.
Let $V$ be a 3 dimensional vectors pace over $\FF_{64}$ equipped with the
natural basis $(e_1,e_2,e_3)$ and
Hermitian form $(x,y) = x_1 y_1^{\tau}+  x_2y_2^{\tau} + x_3y_3^{\tau}$.
Let $\beta=\zeta^7$ and  $\alpha=\beta^{3}$. Then $\beta$ and $\alpha$ 
have order 9 and 3, respectively, and
$\beta\beta^{\tau}=1$ and $\alpha\alpha^{\tau}=1$.
Let $\sigma$ be the automorphism of $\SU_3(8)$ induced by $\rho$.

Let
$$A=\small{\left( \begin{array}{ccc}
0&0&1\\
1&0&0\\
0&1&0\\
\end{array}\right)}, \,\,
B=\small{\left( \begin{array}{ccc}
1&0&0\\
0&\alpha&0\\
0&0&\alpha^{-1}\\
\end{array}\right)}, \,\,
C=\small{\left( \begin{array}{ccc}
\beta&0&0\\
0&\beta^4&0\\
0&0&\beta^4\\
\end{array}\right)}, \,\,
D=\small{\left( \begin{array}{ccc}
1&1&1\\
1&\alpha&\alpha^{-1}\\
1&\alpha^{-1}&\alpha\\
\end{array}\right)},$$

$$
E=\small{\left( \begin{array}{ccc}
1&0&0\\
0&\beta&0\\
0&0&\beta^{-1}\\
\end{array}\right)}, \,\,
F=\small{\left( \begin{array}{ccc}
1&0&0\\
0&0&1\\
0&1&0\\
\end{array}\right)} \,\,
\hbox{ and }Z=\small{\left( \begin{array}{ccc}
\alpha&0&0\\
0&\alpha&0\\
0&0&\alpha\\
\end{array}\right)}.
$$

\par\noindent
Note that $A,B,C,D,E,F,Z \in \SU_3(8)$ with $C^3=Z$, $D^2=F$, $E^3=B$.
For any $X \in {\rm \Gamma U_3(8)}$ we will denote by $\o X$ its image in ${\rm P \Gamma U_3(8)}$.

\par\smallskip\noindent
We define the following groups
$$Q_1=\langle \o A, \o B \rangle, \,\, Q_2= \langle \o A, \o B, \o C\rangle, \,\,
Q_*= \langle \o B, \o C\rangle, \,\, S=\langle \o E, \o F , \o \sigma^3 \rangle,$$
$$H_1=\langle \o A,\o B,\o C,\o D, \o \sigma^3 \rangle , \,\,
H_2=\langle \o A,\o B,\o C,\o E, \o F, \o \sigma^3 \rangle \hbox{ and } H=\langle H_1, H_2 \rangle.$$
\par\noindent
Note that $H \leq \PSU_3(8) \sdp \langle\sigma^3 \rangle$. Moreover let
$$\hat Q_1=\langle \o A, \o B, \o \sigma^2 \rangle, \,\, 
\hat Q_2= \langle \o A, \o B, \o C , \o \sigma^2 \rangle, \,\,
K_1=\langle H_1 ,\o\sigma^2 \rangle , \,\,
K_2=\langle H_2 ,\o\sigma^2 \rangle \hbox{ and } 
K=\langle K_1, K_2 \rangle.$$

\par\noindent
Note that $K \leq \PSU_3(8) \sdp \langle \o\sigma \rangle$ and
$K=\langle H_1, H_2, \o\sigma^2 \rangle=\langle H, \o\sigma^2 \rangle$.

\begin{lemma}\label{L1}
The following relations hold:
\begin{enumerate}[(i)]
\item $[A,B]=Z^2$, $[A,C]=BZ^2$ and $[B,C]=1$;
\item $[D,A]=BA$ and $[D,B]=A^2B$;
\item $A^{\sigma}=A$, $B^{\sigma}=B^{-1}$, $C^{\sigma}=C^{2}$ and $D^{\sigma}=D^{-1}$;
\end{enumerate}
In particular $[\langle \o C, \o D \rangle, Q_1]=Q_1$. 
\end{lemma}

\pf
Straightforward.
\eop

\begin{lemma}\label{L2}
We have
\begin{enumerate}[(i)]
\item $Q_1$ and $Q_*$ are abelian group of order $3^2$ and 
$[Q_1,Q_*]= \langle \o B \rangle$.
\item $Q_2=Q_1Q_*$ is a special group of order $3^3$ and exponent 3 with center $\langle \o B \rangle$;
\item $H_1 \cong P \Gamma U_3(2) \cong \AGL_2(3)$ and $O_3(H_1)=Q_1$;
\item $\hat Q_1=Q_1 \times \langle \o \sigma^2 \rangle$ and 
$\hat Q_2=Q_2 \times \langle \o \sigma^2 \rangle$;
\item $K_1 \cong C_3 \times \AGL_2(3)$ and $O_3(K_1)=\hat Q_1$.
\end{enumerate}
In particular $Z(Q_2) \leq Q_1$ and  $Z(\hat Q_2) \leq \hat Q_1$.
\end{lemma}

\pf
Observe that by \ref{L1}$(i)$ both $Q_1$ and $Q_*$ 
are abelian groups of order $3^2$ and $[Q_1,Q_*]= \langle \o B \rangle$.
Hence $Q_1$ is normalized by $Q_*$ and
$Q_2=Q_1Q_*$. It follows that $Q_2$ is a special group of order $3^3$ and exponent 3.
By \ref{L1} 
$\langle \o C,\o D \rangle$ normalizes $Q_1$.
Since $\o C^4=\o C$ it follows from \ref{L1}$(iii)$ that $\o \sigma^2$
centralizes $\langle \o A,\o B,\o C,\o D \rangle$.
Hence $\langle \o A,\o B,\o C,\o D \rangle\cong PGU_3(2)$.
Moreover \ref{L1}$(iii)$ now also implies
$H_1\cong P \Gamma U_3(2) \cong \AGL_2(3)$. Now $(iv)$, $(v)$ and the last statement
follow easily.
\eop

We define 
$$\Lambda=\{X \leq Q_2 \mid |X|=9, X \hbox{ elementary abelian and } X\neq Q_* \}.$$

Note that, by \ref{L2}, $|\Lambda|=3$ and $Q_1 \in \Lambda$.

\begin{lemma}\label{L3}
The following relations hold:
\begin{enumerate}[(i)]
\item $[E,A]=BC$, $[E,B]=1$ and $[E,C]=1$;
\item $[F,A]=A^2$, $[F,B]=B^2$, $[F,C]=1$  and $[F,E]=E^2$;
\item $E^{\sigma}=E^{2}$ and $F^{\sigma}=F$.
\end{enumerate}
In particular $[\langle \o E \rangle, Q_2]=Q_*$ and 
$[\langle \o E \rangle, \langle \o \sigma^2 , \o B\rangle ]=\langle \o B \rangle$.
\end{lemma}

\pf
Straightforward.
\eop

\begin{lemma}\label{L4}
We have
\begin{enumerate}[(i)]
\item $S= \langle \o E, \o F \rangle \times \langle \o {F\sigma^3} \rangle 
\cong Dih(18) \times C_2$;
\item $S$ normalizes $Q_2$ and $\hat Q_2$, $S \cap Q_2=Z(Q_2)$ 
and $S/Z(Q_2) \cong \Sym(3) \times C_2$.
\item $S$ normalizes $Q_*$, $\langle \o {F\sigma^3} \rangle^{\Lambda}=1$ 
and $S^{\Lambda} \cong \Sym(3)$;
\item $H_2 \cong Q_2\cdot (\AGL_1(3) \times C_2)$ and $H_2/Z(Q_2)\cong \AGL_2(3,S)$;
\item $K_2 \cong \hat Q_2\cdot (\AGL_1(3) \times C_2)$ and 
$K_2/Z(\hat Q_2)\cong \AGL_2(3,S)$.
\end{enumerate}
\end{lemma}

\pf 
$(i)$.
By \ref{L3}$(ii)$ $[F,E]=E^2$, whence $E^F=E^{-1}$ and
$\langle \o E, \o F \rangle \cong Dih(18)$.
It follows from \ref{L3}$(iii)$ that
$\o {F\sigma^3}$ has order 2 and
centralizes $\o E$ and $\o F$.
This implies that 
$S=\langle \o E, \o F \rangle\times \langle \o {F\sigma^3} \rangle 
\cong Dih(18) \times C_2$. 
\par\noindent
$(ii)$. It follows from \ref{L1} and \ref{L3} that $S$
normalizes $Q_2$ and  $[S, \langle \o \sigma^2 \rangle]=\langle \o B \rangle$.
Hence $S$ normalizes $\hat Q_2$.
By \ref{L3}$(i)$ $S$ does not normalize $Q_1$. 
Since $\langle \o E \rangle =O_3(S)$  and
$S \cap Q_2 \leq O_3(S)$ we have $S \cap Q_2 = \langle \o B \rangle =Z(Q_2)$. 
Whence $S/Z(Q_2) \cong \Sym(3) \times C_2$.
\par\noindent
$(iii)$.
By \ref{L3} $S$ normalizes $Q_2$ and $Q_*$. Hence $S$ acts on $\Lambda$.
Since $[\o E, Q_1] \not\leq Q_1$, $\langle \o E \rangle$ is transitive on $\Lambda$.
By $(i)$ $[\o E, \o {F\sigma^3}]=1$, and by \ref{L2}$(iii)$ $F\sigma^3$ normalizes $Q_1$.
It follows that $\langle \o {F\sigma^3} \rangle$ acts trivially on $\Lambda$.
\par\noindent
$(iv)$.
By \ref{L1}$(iii)$ and \ref{L3}$(ii)$
$\o {F\sigma^3}$ inverts $\o A$ and $\o C$, and centralizes $\o B$. 
It follows that
$H_2= Q_2S \cong Q_2 \cdot (\AGL_1(3) \times C_2)$ and $H_2/Z(Q_2)\cong \AGL_2(3,S)$.
\par\noindent
$(v)$. By \ref{L3}$(iii)$ $[Q_2, \langle \o \sigma^2 \rangle ]=1$. Hence
$Z( \hat Q_2)= \langle \sigma^2 , Z(Q_2)\rangle$.
Since $S$ normalizes $\hat Q_2$ it follows that
$K_2 =  \hat Q_2 S \cong \hat Q_2\cdot (\AGL_1(3) \times C_2)$ and 
$K_2/Z( \hat Q_2)\cong \AGL_2(3,S)$.
\eop

\begin{lemma}\label{L5} We have
\begin{enumerate}[(i)]
\item $H_1 \cap H_2 = \langle \o A,\o B,\o C,\o F, \o \sigma^3 \rangle
\cong \AGL_2(3,S)$ and $H \cong \PSU_3(q) \sdp \langle \o \sigma^3 \rangle$;
\item $K_1 \cap K_2 = \langle \o A,\o B,\o C,\o F, \o \sigma^3, \o \sigma^2 \rangle
\cong C_3 \times \AGL_2(3,S)$ and $K \cong \PSU_3(q) \sdp \langle \o \sigma\rangle$.
\end{enumerate}
\end{lemma} 

\pf
$(i)$. Let $J= \langle \o A,\o B,\o C,\o F, \o \sigma^3 \rangle$.
Clearly $J \leq H_1 \cap H_2$ and $J \cong Q_2 \sdp (C_2 \times C_2)$.
Since $\o D$ does not normalize $Q_2$ we have that
$H_1 \cap H_2$ is a proper subgroup of $H_1$.
Since $J/Q_1 \cong C_2 \times \Sym(3)$ is a maximal subgroup of 
$H_1 /Q_1 \cong GL_2(3)$  it follows that $J=H_1 \cap H_2 \cong AGL_2(3,S)$.
By \ref{L2}$(iii)$ and \ref{L4}$(iv)$
$H_1 \cong P\Gamma U_3(2)$, $O_3(H_1)=Q_1$ and
$H_2 \cong  Q \cdot (C_2 \times \Sym(3))$, where $Q_2$ a special 3 group of order $3^3$.
Since $H_1$ is a maximal subgroup of $\PSU_3(8) \sdp \langle \o \sigma^3 \rangle$,
see \cite{atlas}, and $H_2 \not \leq H_1$ we
have $H\cong \PSU_3(8) \sdp \langle \o \sigma^3 \rangle$. 
\par\noindent
$(ii)$.
By \ref{L1} and \ref{L3} 
$\o \sigma^2$ centralizes $H_1$ and normalizes $H_2$. 
Hence $\o \sigma^2$ normalizes $H$ and thus
$K\cong \PSU_3(8) \sdp \langle \o \sigma \rangle$.
Moreover $|K_i:H_i|=3$ and 
$\langle \o \sigma^2 , J \rangle\leq K_1 \cap K_2$.
It follows that 
$K_1 \cap K_2 = \langle \o \sigma^2 , J \rangle \cong C_3 \times AGL_2(3,S)$.
\eop

\begin{lemma} \label{L6} 
The following hold:
\begin{enumerate}[(i)]
\item The amalgam $H_1 \leftarrow H_1 \cap H_2 \rightarrow H_2$ is of shape ${\cal D}_2$;
\item The amalgam $K_1 \leftarrow K_1 \cap K_2 \rightarrow K_2$ is of shape ${\cal E}_2$.
\end{enumerate}
\end{lemma}

\pf
Recall the notation of section 2.
\par\noindent
$(i)$.
Let $H_{1,2}=H_1 \cap H_2$ and let 
${\mathfrak A}: H_1 \leftarrow H_{1,2} \rightarrow H_2$.
It follows from \ref{L2}, \ref{L3}, \ref{L4}$(iv)$ and \ref{L5}$(i)$ that
$T_1=\langle \o A,\o B,\o F \rangle$ and
$T_2=\langle \o A,\o B,\o C,\o {F\sigma^3} \rangle$.
Note that 
$\langle T_1,T_2 \rangle = H_{1,2}$ and thus $X=H_{1,2}$.
Now $O_3(X)=Q_2$ and $Z(Q_2)=\langle \o B \rangle$. By  \ref{L3} we have
that $C_X(Z(O_3(X)))=\langle \o A,\o B,\o C,\o {F\sigma^3} \rangle$. 
Hence $C_X(Z(O_3(X)))=T_2$.
By \ref{L4}$(iv)$ $H_2/Z(Q_2)\cong \AGL_2(3,S)$. Whence
the amalgam ${\mathfrak A}$ has shape $\AGL_2(3,S)$.
We have $H_1 \cong \AGL_2(3)$ and $H_{1,2}\cong\AGL_2(3,S)$ by
\ref{L2}$(iii)$ and \ref{L5}$(i)$, respectively.
By \ref{L4}$(iv)$ $H_2 \cong Q_2\cdot (\AGL_1(3) \times C_2)$.
Note that $T_2 \cong \AGL_2(3,S)^{\#}$, 
$O_3(H_2)=\langle \o A,\o B,\o C, \o E \rangle$ and, by \ref{L3},
$C_{O_3(H_2)}(\o {F\sigma^3})= \langle \o E\rangle$. 
It follows  that amalgam ${\mathfrak A}$ has shape ${\cal D}_2$.
\par\noindent
$(ii)$.
Let $K_{1,2}=K_1 \cap K_2$ and let ${\mathfrak B}: K_1 \leftarrow K_{1,2} \rightarrow K_2$.
It follows from \ref{L2}, \ref{L3}, \ref{L4}$(v)$ and \ref{L5}$(ii)$ that
$T_1=\langle  \o \sigma^2, \o A,\o B,\o F \rangle$ and
$T_2=\langle  \o \sigma^2, \o A,\o B,\o C,\o {F\sigma^3} \rangle$.
Hence $\langle T_1,T_2 \rangle = K_{1,2}$ and thus $X=K_{1,2}$.
Now $O_3(X)=\hat Q_2$ and $Z(O_3(X))=Z(\hat Q_2)=\langle \o \sigma^2, \o B \rangle$.
By  \ref{L3} we have
that $C_X(Z(O_3(X)))=\langle \o \sigma^2,\o A,\o B,\o C,\o {F\sigma^3} \rangle$.
Hence $C_X(Z(O_3(X)))=T_2$.
By \ref{L4}$(v)$ $K_2 /Z(\hat Q_2) \cong \AGL_2(3,S)$.
Whence the amalgam ${\mathfrak B}$ has shape $\AGL_2(3,S)$.
We have $K_1 \cong \AGL_2(3)$ and $K_{1,2} \cong C_3 \times \AGL_2(3,S)$ by
\ref{L2}$(v)$ and \ref{L5}$(ii)$, respectively.
By \ref{L4}$(v)$ $K_2 \cong \hat Q_2\cdot (\AGL_1(3) \times C_2)$.
Note that $Z(T_2)= \langle \o \sigma^2, \o B \rangle$ and
$T_2 \cong C_3 \times \AGL_2(3,S)^{\#}$. Since
$K_2/C_X(Z(O_3(X))) \cong \langle \o E, \o F \rangle/  \langle \o B \rangle$ we get
$Z(O_3(X)) \sdp K_2/C_{K_2}(Z(O_3(X))) \cong \AGL_2(3,S)^*$.
Moreover,
$O_3(K_2)=\langle \o \sigma^2, \o A,\o B,\o C, \o E \rangle$.
Now \ref{L3} gives
$C_{O_3(K_2)}(\o {F\sigma^3})= \langle \o \sigma^2, \o E\rangle \cong SP_2$. 
It follows that the amalgam ${\mathfrak B}$  has shape ${\cal E}_2$.
\eop

Let $\Delta$ be the coset graph on
$K_1$ and $K_2$ with the edge stabilizer $K_1\cap K_2$. Then
$\Delta$ is a biregular bipartite graph with valencies $4$ and $3$.
The group $K$ is transitive on the two sets of vertices of equal valency, and $\Delta$
is locally $s$-arc transitive with $s \geq 1$. By \ref{L6} the vertex stabilizer
amalgam will be of shape ${\cal E}_2$. Now \cite[3.2]{PushingUpType1} implies that the resulting graph will have $s=5$. We will give an alternative proof of this fact here
to make the paper self contained and to obtain additional information on the stabilizer
of a 5-arc.

\begin{prop}\label{L7}
$\Delta$ is a locally 5-arc transitive $K$-graph of pushing up type of valencies 3 and 4, 
where $K \cong \PSU_3(q) \sdp \langle \o \sigma\rangle$. Moreover, for any
$5$-arc $\alpha = (y_0,y_1,y_2,y_3,y_4,y_5)$ we have $K_{\alpha}=Z(O_3(K_{y_2,y_3}))$.
\end{prop}

\pf
By \ref{L5}$(ii)$ $K \cong \PSU_3(q) \sdp \langle \o \sigma\rangle$.
Let $x_1=K_1$ and $x_2=K_2$ be an edge.
Then $O_3(K_1) \leq Q_3(K_2)$ and $C_{K_i}(O_3(K_i)) \leq O_3(K_i)$ for $i\in\{1,2\}$.
Hence $\Delta$ is of pushing up type.

We define
$x_3=K_1 \o E$, $x_0=K_2  \o D$, $x_4=K_2  \o {DE}$ and $x_{-1}=K_1 \o {ED}$.
Hence $x_{-1}, x_0,x_1,x_2,x_3,x_4$ is a 5-arc.
Then
$$K_{x_3}=(K_1)^{\o E}, \,\, K_{x_0}=(K_2)^{\o D}, \,\, K_{x_4}=(K_0)^{\o E}, \hbox{ and } K_{x_{-1}}=(K_3)^{\o D}.$$

\par\smallskip\noindent
By \ref{L2}$(v)$, \ref{L4}$(v)$ and \ref{L5}$(ii)$ we have
$|K_{x_1}:K_{x_1,x_2}|=4$ and $|K_{x_2}:K_{x_1,x_2}|=3$. Whence $s \geq 1$.
\par\smallskip\noindent
Since $\hat Q_2 \langle \o {F\sigma^3}\rangle \leq K_{x_1,x_2}$ is normalized by
$\o E$ we have $\hat Q_2 \langle \o {F\sigma^3}\rangle \leq K_{x_2,x_3}$.
Observe that $Z(\hat Q_2) \leq \hat Q_1=O_3(K_{x_1}) \leq \hat Q_2$.
Since $\langle \o B , \o \sigma^2 \rangle=Z(\hat Q_2) \nl K_{x_2}$ it follows
from \ref{L4} that the action of $K_{x_2}$ on $\Delta (x_2)$ is 
equal to that of $S$ on $\Lambda$.
Hence we have that
$K_{x_1,x_2} \cap K_{x_2,x_3} \leq \hat Q_2 \langle \o {F\sigma^3}\rangle$.
It follows that
$$K_{x_1,x_2,x_3}=K_{x_1,x_2} \cap K_{x_2,x_3} =
\langle \o A,\o B,\o C, \o {F\sigma^3}, \o \sigma^2 \rangle 
\cong \hat Q_2 \sdp \langle \o {F\sigma^3} \rangle .$$
\par\noindent
Since $K_{x_0,x_1}=K_{x_1,x_2}^{\o D}= 
\langle \o A,\o B,\o C^{\o D},\o F,\o {F\sigma^3}, \o \sigma^2  \rangle$
we have
$\hat Q_1 \langle \o F, \o \sigma^3 \rangle \leq K_{x_0,x_1} \cap K_{x_1,x_2}$. 
Since $K_1/{\hat Q_1} \cong GL_2(3)$ and $\o D$ does
not normalize $\hat Q_2 /{\hat Q_1}$ we have
$$K_{x_0,x_1,x_2}=K_{x_0,x_1} \cap K_{x_1,x_2}=
\langle \o A,\o B,\o F,\o \sigma^3, \o \sigma^2 \rangle 
\cong \hat Q_1 \sdp \langle \o F,\o \sigma^3 \rangle.$$
\par\noindent
In particular $|K_{x_1,x_2}:K_{x_1,x_2,x_3}|=2$ and 
$|K_{x_1,x_2}:K_{x_0,x_1,x_2}|=3$. 
\par\noindent
Since $s \geq 1$ it follows that for $i \in\{1,2\}$,
$G_{x_i}$ is transitive on 2-arcs originating at $x_i$.
Whence $s \geq 2$.  
\par\smallskip\noindent
Since $\hat Q_1 \langle \o {F\sigma^3}\rangle 
\leq K_{x_0,x_1,x_2} \cap K_{x_1,x_2,x_3}$, 
$|K_{0,1,2}|=3^3\cdot 2^2$ and $|K_{x_1,x_2,x_3}|=3^4 \cdot 2$  
it follows  that
$$K_{x_0,x_1,x_2,x_3}=K_{x_0,x_1,x_2} \cap K_{x_1,x_2,x_3} =
\langle \o A,\o B, \o \sigma^2, \o {F\sigma^3}  \rangle
\cong \hat Q_1 \sdp \langle \o {F\sigma^3}\rangle.$$
\par\noindent
In particular $|K_{x_0,x_1,x_2}:K_{x_0,x_1,x_2,x_3}|=2$ and 
$|K_{x_1,x_2,x_3}:K_{x_0,x_1,x_2,x_3}|=3$. 
\par\noindent
Since $s \geq 2$ it follows that for $i \in\{0,3\}$,
$G_{x_i}$ is transitive on 3-arcs originating at $x_i$.
Whence $s \geq 3$.  
\par\smallskip\noindent
Observe that $\o E$ normalizes $K_{x_1,x_2,x_3}$ and that
$K_{x_0,x_1,x_2,x_3,x_4} =  K_{x_0,x_1,x_2} \cap K_{x_1,x_2,x_3} \cap K_{x_2,x_3,x_4}$.
Since $K_{x_2,x_3,x_4}=K_{x_0,x_1,x_2}^{\o E}$,
$K_{x_0,x_1,x_2} \cap K_{x_1,x_2,x_3}= \hat Q_1 \langle \o {F\sigma^3}\rangle$ 
and $\hat Q_1 \cap \hat Q_1^{\o E}=\langle \o B, \o \sigma^2 \rangle$, 
it follows  that
$$K_{x_0,x_1,x_2,x_3,x_4}=K_{x_0,x_1,x_2}\cap K_{x_2,x_3,x_4}= \langle \o B, \o \sigma^2, \o {F\sigma^3}\rangle 
\cong \langle \o B, \o \sigma^2 \rangle \times \langle \o {F\sigma^3}\rangle.$$
Since
$K_{x_{-1},x_0,x_1}=K_{x_1,x_2,x_3}^{\o D} =
\langle \o A,\o B,\o C^{\o D}, \o \sigma^3 , \o \sigma^2 \rangle$ it follows that
$\hat Q_1 \leq K_{x_{-1},x_0,x_1}\cap K_{x_0,x_1,x_2,x_3}$. 
Since $K_{x_1}/\hat Q_1\cong GL_2(3)$ and 
$K_{x_0,x_1,x_2,x_3}=\hat Q_1 \sdp \langle \o{F\sigma^3}\rangle$ it follows that
$$K_{x_{-1},x_0,x_1,x_2,x_3}=K_{x_{-1},x_0,x_1}\cap K_{x_0,x_1,x_2,x_3}= 
\langle \o A, \o B, \o \sigma^2 \rangle= \hat Q_1.$$
\par\noindent
In particular $|K_{x_0,x_1,x_2,x_3}:K_{x_{-1},x_0,x_1,x_2,x_3}|=2$ and 
$|K_{x_1,x_2,x_3,x_4}:K_{x_0,x_1,x_2,x_3,x_4}|=3$. 
\par\noindent
Since $s \geq 3$ it follows that for $i \in\{3,4\}$,
$G_{x_i}$ is transitive on 4-arcs originating at $x_i$.
Whence $s \geq 4$.  
\par\smallskip\noindent
Moreover
$$K_{x_{-1},x_0,x_1,x_2,x_3,x_4}=K_{x_{-1},x_0,x_1,x_2,x_3}\cap K_{x_0,x_1,x_2,x_3,x_4}= \langle \o B, \o \sigma^2 \rangle.$$
\par\noindent
In particular $|K_{x_0,x_1,x_2,x_3,x_4}:K_{x_{-1},x_0,x_1,x_2,x_3,x_4}|=2$ and 
$|K_{x_{-1},x_0,x_1,x_2,x_3}:K_{x_{-1},x_0,x_1,x_2,x_3,x_4}|=3$. 
\par\smallskip\noindent
Since $s \geq 4$ it follows that for $i \in\{4,-1\}$,
$G_{x_i}$ is transitive on 5-arcs originating at $x_i$.
Whence $s \geq 5$.  
\par\smallskip\noindent
Since $K_{x_{-1},x_0,x_1,x_2,x_3,x_4}$ has odd order $s \leq 5$.
Moreover $K_{x_{-1},x_0,x_1,x_2,x_3,x_4}= Z(\hat Q_2)=Z(O_3(K_{x_1,x_2}))$.
Since any 5-arc is conjugated to either $x_{-1},x_0,x_1,x_2,x_3,x_4$ or $x_4,x_3,x_2,x_1,x_0,x_{-1}$ the
proposition follows.
\eop

\par\smallskip
Next we give the structure of the vertex stabilizers in terms of the kernels of action.

\begin{lemma}\label{L8}
$\Delta$ is of local characteristic $3$ with respect to $K$.
Moreover, the following hold:
\begin{enumerate}[(i)]
\item $K_{x_1}/K_{x_1}^{[1]} \cong \Sym(4)$, $K_{x_1}^{[1]} \cong \hat Q_1 \sdp C_2$, $K_{x_1}^{[2]} \cong \hat Q_1$, $K_{x_1}^{[3]} =K_{x_1}^{[4]}=Z(K_{x_1})\cong C_3$, $K_{x_1}^{[5]}=1$;
\item $K_{x_2}/K_{x_2}^{[1]} \cong \Sym(3)$, $K_{x_2}^{[1]} \cong \hat Q_2 \sdp C_2$, $K_{x_2}^{[2]}=K_{x_2}^{[3]} \cong Z(\hat Q_2) \cong C_3 \times C_3$, $K_{x_2}^{[4]}=1$.
\end{enumerate}
\end{lemma}
 
\pf
Let $\{x_1,x_2\}$ be an edge. Since $K$ is edge transitive we may assume that $K_{x_1}=K_1$ and $K_{x_2}=K_2$.
Since $\Delta$ is biregular of valency 3 and 4 it follows from
\ref{L2}$(v)$, \ref{L4}$(v)$ and \ref{L5}$(ii)$ that $K_{x_1}/K_{x_1}^{[1]} \cong \Sym(4)$,
$K_{x_2}/K_{x_2}^{[1]} \cong \Sym(3)$, $K_{x_1}^{[1]} ={\hat Q_1} \sdp \langle \o {F}\rangle$ and
$K_{x_2}^{[1]} = {\hat Q_2} \sdp \langle \o {F\sigma^3}\rangle$.
Consequently $K_{x_1,x_2}^{[1]}={\hat Q_1}=K_{x_1}^{[2]}$.
Since $K_{x_1,x_2}^{[1]}=K_{x_1}^{[2]}$ we have
$K_{x_1}^{[3]}=K_{x_1}^{[4]}$ and $K_{x_2}^{[2]}=K_{x_2}^{[3]}$.
Moreover \ref{L2} yields
$C_{K_1}(\hat Q_1) \leq C_{K_1}(Q_1) = \hat Q_1$.
By \ref{L2}$(ii)$ and \ref{L2}$(iv)$ $Q_2$ is an extraspecial group and
$Z(\hat(Q_2))=Z(Q_2) \times \langle \o \sigma^2 \rangle$.
Now \ref{L4}$(iv)$ yields
$C_{K_2}(\hat Q_2) \leq C_{K_2}(Q_2) = \hat Q_2$. Hence
$\Delta$ is of local characteristic 3, and
of pushing up type since ${\hat Q_1}\leq {\hat Q_2}$.
In particular we have $Z({\hat Q_2})\leq {\hat Q_1}$. 
Whence $Z({\hat Q_2})\leq K_{x_2}^{[2]}$. 
Suppose $K_{x_2}^{[2]}={\hat Q_2}$. Then ${\hat Q_2}\leq O_p(K_{x_1}^{[1]})=Q_1$
which is a contradiction.
Hence $K_{x_2}^{[2]} \neq {\hat Q_2}$.
By \ref{L2}$(ii)$ and \ref{L2}$(v)$ ${\hat Q_2}={\hat Q_1}Q_*$ and it now follows from
\ref{L4}$(v)$ that $Q_*Z({\hat Q_2})/Z({\hat Q_2})$
is the only proper non-trivial $S$ invariant subgroup of ${\hat Q_2}/Z({\hat Q_2})$.
Since $Q_* \not\leq K_{x_1}^{[1]}$ it follows that $Z({\hat Q_2})= K_{x_2}^{[2]}$.
Since $K_{x_1}^{[3]} \leq  K_{x_2}^{[2]}$ is normalized by $K_{x_1}$, we have
$K_{x_1}^{[3]}= \langle \o {\sigma^2}\rangle$. Hence $K_{x_1}^{[3]}= Z(K_{x_1})$.
Since $K_{x_2}^{[4]}$ is an $S$-invariant subgroup of $K_{x_1}^{[3]}$ we must have
$K_{x_2}^{[4]}=1$, and consequently $K_{x_1}^{[5]}=1$ too.
\eop

\begin{lemma}\label{L9}
Let  $\alpha=(y_0,y_1,y_2,y_3,y_4,y_5)$ be a 5-arc with $|\Delta (y_0)|=3$.
Then $K_{\alpha}$ is transitive on $\Delta (y_5) \backslash \{y_4\}$.
In particular, if $y$ is a vertex of $\Delta$ with $|\Delta(y)|=3$, then $K_{y}$ is transitive on 6-arcs originating at $y$.
\end{lemma}

\pf
By \ref{L7} and \ref{L5}$(ii)$ $K_{\alpha}=Z(O_3(K_{y_2,y_3})) \cong Z(\hat Q_2)\cong C_3 \times C_3$.
Since $|\Delta(y_0)|=3$ we have $Z(O_3(K_{y_2,y_3}))=Z(O_3(K_{y_2}^{[1]}))$.
Moreover, since $s \geq 1$, $Z(O_3(K_{y_2,y_3}))$ is not normalized by $K_{y_3}$.
In particular, $Z(O_3(K_{y_2}^{[1]}))\neq Z(O_3(K_{y_4}^{[1]}))$.
By
\ref{L2} $Z(O_3(K_{y_4}^{[1]})) \leq O_3(K_{y_3}^{[1]})$ and $|O_3(K_{y_3}^{[1]}):Z(O_3(K_{y_4}^{[1]})) |=3$ 
it follows that
$O_3(K_{y_3}^{[1]})=\langle Z(O_3(K_{y_2}^{[1]})), Z(O_3(K_{y_4}^{[1]})) \rangle$.

Suppose $Z(O_3(K_{y_2}^{[1]})) \leq K_{y_5} ^{[1]}$. Then $Z(O_3(K_{y_2}^{[1]})) \leq O_3( K_{y_5}^{[1]})$,
and thus we have 
$$\langle Z(O_3(K_{y_2}^{[1]})), Z(O_3(K_{y_4}^{[1]})) \rangle \leq O_3( K_{y_5}^{[1]}).$$
Hence $O_3( K_{y_3}^{[1]})=O_3(K_{y_5}^{[1]})$. 
Since $K_{y_4}= \langle K_{y_3,y_4}, K_{y_4,y_5}\rangle $ we get 
$O_3( K_{y_3}^{[1]}) \nl K_{y_4}$. 
Hence $O_3( K_{y_3}^{[1]})$
is normalized by $\langle  K_{y_3}, K_{y_4}\rangle$.
It follows that  $O_3( K_{y_3})=1$. A contradiction.
Hence $Z(O_3(K_{y_2}^{[1]}))\not\leq K_{y_5} ^{[1]}$ and thus induces a group of order 3 on
$\Delta(y_5) \backslash \{y_4\}$. In particular $K_{\alpha}$ is transitive on
$\Delta(y_5) \backslash \{y_4\}$, and the lemma follows since $s=5$.
\eop

We finish our investigation of the $K$-graph $\Delta$ by showing that $K$ is actually the full automorphism group of $\Delta$.

\begin{lemma}\label{L10}
$\Aut (\Delta) =K$.
\end{lemma}

\pf Let $G=\Aut (\Delta)$. Since $\Delta$ is biregular
with valencies 3 and 4 it follows that $G_z$ is a $2,3$-group for any $z \in V\Delta$.
Since $1 \neq K_z^{[2]} \leq G_z^{[2]}$
it follows from \cite[3.2]{NonLocalChar} that $\Delta$ is of local characteristic $p$, and by
\cite{TW} that $p=3$. Let $\{x_1,x_2\}\in E\Delta$, with $|\Delta(x_1)|=4$.
Since $|\Delta (x_2)|=3$ it follows 
from \cite[2.6 and 2.7]{TW} that $G_{x_1,x_2}^{[1]} =G_{x_1}^{[2]}$ and
$G_{x_2}^{[2]}=G_{x_2}^{[3]}$ is a $3$-group.

For  $z \in V\Delta$ we define $Q_z=O_3(G_z^{[1]})$.
Let $L_{x_1} = \langle Q_u \mid u \in \Delta (x_1) \rangle Q_{x_1}$.
Note that
$Q_{x_1} \leq Q_{x_2}$, so $Q_x \leq G_{x_1} ^{[2]}$, and $L_{x_1} ^{\Delta (x_1)} \cong PSL_2(3)$.
Let $x_0 \in \Delta (x_1) \setminus \{x_2\}$.
Then $L_{x_1} \cap G_{x_0,x_1,x_2} = L_{x_1}\cap G_{x_1}^{[1]}$. 
Since $s \geq 4 $ it follows from \cite[3.5]{NonLocalChar} that
$|(L_{x_1}\cap G_{x_1}^{[1]})/Q_{x_1}| \neq 1$.
Let $N=O^3(L_{x_1})$. Then $|N^{\Delta (x_1)}|=4$ and it follows from
\cite[5.2]{bonstell} that $NQ_{x_1}/Q_{x_1} \cong Q_8$.

Let $t\in N$ be an involution. 
Then $L_{x_1}\cap G_{x_1}^{[1]}=\langle t \rangle Q_{x_1}$, hence
$t \not \in G_{x_1,x_2}^{[1]}$.
Let  $L_{x_2}= \langle t^{G_{x_2}} \rangle Q_{x_2}$. Then
$L_{x_2}^{\Delta (x_2)} \cong \Sym(3)$. Since
$\langle t \rangle Q_{x_1} \nl G_{x_1}$ it follows that
$[t, G_{x_2}^{[1]}] \leq Q_{x_1} \leq Q_{x_2}$.
Hence $[L_{x_2},G_{x_2}^{[1]}] \leq Q_{x_2}$.
By \cite[5.1]{bonstell} $[L_{x_1},G_{x_1}^{[1]}] \leq Q_{x_1}$.
Let $S_i \in Syl_2 (G_{x_i}^{[1]})$, $i=1,2$. By the Frattini argument
$N_{L_{x_i}}(S_i)$ is transitive on $\Delta (x_i)$.
Since $[N_{L_{x_i}}(S_i),S_i] \leq S_i \cap Q_{x_i}=1$ we have
$C_{L_{x_i}}(S_i)$ is transitive on $\Delta (x_i)$ for $i=1,2$. Hence
$S_1 \cap S_2=1$ and so $G_{x_1,x_2}^{[1]}$ is  $3$-group. It follows that 
$G_{x_2}/Q_{x_2} \cong C_2 \times \Sym(3)$, hence $G_{x_2}$ is the normalizes
a Sylow $3$-subgroup of $G$. Moreover
Sylow $2$-subgroups of $G_{x_1}$ and $G_{x_2}$ are of order 16 and 4, respectively.

Note that $\Delta$ is bipartite with blocks $\Delta_1$ and $\Delta_2$, where
$|\Delta_1|= 2^6\cdot 3 \cdot 7 \cdot 19$ and $|\Delta_2|= 2^8 \cdot 7 \cdot 19$.
In particular a Sylow $2$-subgroups of $G$ has order $2^{10}$.
Let $z \in \Delta_2$. Since $K_z$ is a maximal subgroup of $K$, 
$K$ is primitive on $\Delta_2$. Whence
$G$ is primitive on $\Delta_2$ too.
Since $|\Delta_2|$ is not a power of a smaller number and since 
in a group of order $|\Delta_2|$ any Sylow $19$-subgroup is normal,
hence is not simple, it follows from the O'Nan-Scott Theorem, 
see for example \cite{LPS1}, that $G$ is almost simple. 
Since $5$ does not divide $|\Delta_2|$, it follows that $5$ does not
divide $|G|$.
Let $L=Soc (G)$. Then $L$ is of Lie-type in characteristic $p$, with $p\neq 5$.
Since, for any prime $p$ different from $5$, 
$5$ divides $q^4-1$, we have
$$L \in \{PSL_2(q), PSL_3(q), \PSU_3(q), G_2(q), {}^2B_2(q), {}^3D_4(q), {}^2G_2(q)\}.$$

Since $\Aut (L)/L$ is solvable we have $Soc(K)\leq L$ and hence
$L$ is primitive on $\Delta_2$.
It follows that $G_z\cap L$ is a maximal subgroup of $L$ 
and the normalizer of a Sylow $3$-subgroup of $L$. 
If $p = 3$, then $G_z \cap L$ is a maximal parabolic subgroup of $L$ and it
follows that $q=3$.  A contradiction with $Soc(K)\leq L$.
Whence $p \neq 3$. 
Since $|G|=2^{10} \cdot 7 \cdot 19 \cdot |S|$, where
$S \in Syl_3(G)$ this yields $p=2$ and $L \cong \PSU_3(8)$. It follows from the
maximal subgroup structure of $\Aut(\PSU_3(8))$ that $G=K$. 
\eop

The following proposition shows that the action of the subgroup $H$ on $\Delta$ is locally 5-arc transitive too.

\begin{prop}\label{L11}
$\Delta$ is a locally 5-arc transitive $H$-graph of pushing up type, where $H \cong \PSU_3(q) \sdp \langle \o \sigma^3 \rangle$.
Moreover the following hold:
\begin{enumerate}[(i)]
\item Let $\{x_1,x_2\}$ be an edge of $\Delta$. Then
\begin{enumerate}[(1)]
\item $H_{x_1}/H_{x_1}^{[1]} \cong \Sym(4)$, $H_{x_1}^{[1]} \cong Q_1 \sdp C_2$, $H_{x_1}^{[2]} \cong Q_1$, $H_{x_1}^{[3]}=1$;
\item $H_{x_2}/H_{x_2}^{[1]} \cong \Sym(3)$, $H_{x_2}^{[1]} \cong Q_2 \sdp C_2$, $H_{x_2}^{[2]}=H_{x_2}^{[3]} =Z(Q_2)\cong C_3$, $H_{x_2}^{[4]}=1$.
\end{enumerate}
\item Let $y$ be a vertex with $|\Delta(y)|=4$, then $H_{y}$ is transitive on 6-arcs originating at $y$;
\end{enumerate}
\end{prop}

\pf
By \ref{L5}$(ii)$ $H \cong \PSU_3(q) \sdp \langle \o \sigma^3 \rangle$.
Let $\{x_1,x_2\}$ be an edge such that $K_{x_1}=K_1$ and $K_{x_2}=K_2$.
Then $H_{x_i}=K_i \cap H$.
Observe that $\o \sigma^2 \not \in H$ and normalizes $H$. Thus we have $H \nl K$, $|K:H|=3$
and $H_i \nl K_i$ with $|K_i:H_i|=3$ for $i=1,2$. 
We conclude that $H_i=K_i \cap H=H_{x_i}$.
For $i \in \{1,2\}$ and $j\in \{1,2,3,4\}$ we have $H_{x_i}^{[j]}=K_{x_i}^{[j]} \cap H$.
Now $(i)$ follows from \ref{L8}, since $\o \sigma^2 \not\in H$.

By \ref{L2}$(iii)$ $C_{H_1}(Q_1) = Q_1$, and by \ref{L4}$(iv)$ 
$C_{H_2}(Q_2) \leq Q_2$.
It follows that
$\Delta$ is a $H$-graph of local characteristic 3, and 
since $O_3(H_{x_1}) \leq O_3(K_{x_1})$,
$\Delta$ is a $H$-graph of pushing up type.

Let $\alpha= (x_{-1},x_0,x_1,x_2,x_3,x_4)$ be a 5-arc containing $(x_1,x_2)$ in its center.
By \ref{L7} $K_{\alpha}=Z(\hat Q_2)= \langle \o B , \o \sigma^2 \rangle$.
Since $\langle \o B \rangle \leq K_{\alpha} \cap H$  
it follows that $K_{\alpha} \cap H= \langle \o B \rangle$.
Hence $H_{\alpha}= \langle \o B \rangle = Z(O_3(H_{x_1,x_2}))$.
In particular $|K_{\alpha}:H_{\alpha}|=3$.
Hence $|K_{x_{-1}}:K_{\alpha}|=|H_{x_{-1}}:H_{\alpha}|$ and 
$|K_{x_{4}}:K_{\alpha}|=|H_{x_{4}}:H_{\alpha}|$.
It follows that $H$ is transitive on the set of 5-arcs with initial vertex $x_{-1}$
and the set of 5-arcs with initial vertex $x_4$. 
Hence $\Delta$ is a locally 5-arc transitive $H$-graph.

By \ref{L8}$(i)$ $K_{x_1} ^{[4]}= Z(K_{x_1}) = \langle \o \sigma^2 \rangle$,
and by  \ref{L7} $K_{\alpha}=Z(O_3(K_{x_1,x_2}))=\langle \o B, \sigma^3\rangle$.
Hence by \ref{L9} $K_{\alpha} \cap K_{x_{-2}}= \langle \o \sigma^2 \rangle$.
It follows that $H_{\alpha} \cap H_{x_{-2}}= 1$, 
$H_{\alpha}$ is transitive on $\Delta (x_{-1}) \backslash \{x_0\}$ and thus $(ii)$ holds. 
\eop

\medskip\noindent
{\it Proof of the main theorems}.
\par
By construction $\Delta$ is a bipartite biregular graph with valencies 3 and 4.
Theorem \ref{main1} follows from \ref{L6}, \ref{L7}, 
 \ref{L9}, \ref{L10} and \ref{L11}$(ii)$.
In order to prove  Theorem \ref{main2} first observe that \ref{L2} and \ref{L4} 
give the structure of $W_i$ and $\hat W_i$, $i \in \{1,2\}$. Now \ref{L8} 
and \ref{L11}$(i)$ complete the proof.

\par\medskip\noindent
{\bf Disclosure statement}
\par\smallskip\noindent
There are no competing interests to declare.

\vspace{1cm}

\end{document}